\documentstyle[11pt,amstex,amssymb]{amsart}\textheight=19.8truecm\textwidth=14.1truecm
\begin{document}

 \centerline{\huge \bf A matrix subadditivity inequality}
\centerline{\huge \bf for symmetric norms}

\vskip 25pt \centerline{{\large \sl \ Jean-Christophe Bourin }}

\vskip 15pt\centerline{\it Dedicated to Fran\c{c}oise Lust-Piquard, with affection}

\vskip 30pt \noindent {\small {\bf Abstract.}  Let $f(t)$ be a
non-negative concave function on $[0,\infty)$.
 We prove that
$$
\Vert\, f(|A+B|)\,\Vert \le \Vert\, f(|A|)+f(|B|)\,\Vert
$$
for all normal $n$-by-$n$ matrices $A$,  $B$  and all symmetric
norms. This  result has several applications. For instance, for a
Hermitian ${\mathbb{A}}=[A_{i,\,j}]$  partitioned in blocks of same
size,
$$
\left\|\, f(|{\mathbb{A}}|) \,\right\| \le \left\|\,\sum
f(|A_{i,\,j}|) \,\right\|
$$
We also prove, in a similar way, that given  $Z$  expansive  and
 $A$ normal of same size,
$$\Vert\, f(|Z^*AZ|)\,\Vert \le \Vert\, Z^*f(|A|)Z\,\Vert.$$


Mathematical subjects classification:   15A60, 47A30, 47A60}

\vskip 30pt\noindent {\Large\bf  1. Some recent results for positive
operators} \vskip 10pt

Several nice   inequalities for concave functions of operators have
been recently established in a serie of papers \cite{B2}, \cite{BU},
\cite{BL} and \cite{B3}. Most of these results are matrix versions
of the obvious
 inequality
 \begin{equation}
 f(a+b) \le f(a)+f(b)
 \end{equation}
 for  non-negative concave functions $f$ on $[0,\infty)$  and
 scalars $a,\,b\ge0$.  By matrix version we mean suitable extension where scalars are replaced by $n$-by-$n$
 matrices, i.e., operators on an
  $n$-dimensional Hilbert space ${\mathcal{H}}$.  For instance, we have \cite{BU}:

\vskip 10pt\noindent
 {\bf Theorem 1.1.} {\it Let $A,\,B\ge0$ and let
 $f:[0,\infty)\to [0,\infty)$ be concave. Then, for all symmetric norms,
 $$\Vert\, f(A+B)\,\Vert \le \Vert\, f(A)+f(B)\,\Vert.$$
 }

\vskip 5pt\noindent As usual, capital letters $A,\,B,\dots $ stand
for operators, $A\ge 0$ refers to positive semi-definite, and a
symmetric norm
  (or unitarily invariant)
  satisfies  $\Vert A\Vert = \Vert UAV\Vert$ for all $A$ and all
unitaries $U,\, V$. Thus,  up to symmetric norms, the basic
inequality (1) still holds on the cone of positive operator. This
subadditivity result for norms can not be extended to the
determinant, even in the case of an operator concave function such
as $f(t)=\sqrt{t}$. The most elementary case in the above theorem is
for the trace norm. Then, the result can be restated as a famous
trace inequality \cite{Ro}:

\vskip 10pt\noindent {\bf Rotfel'd Inequality.} {\it Let $f$ be a
concave function on $[0,\infty)$ such that $f(0)\ge0$. Then, for all
$A,\,B\ge0$,
$${\mathrm{Tr}}\, f(A+B)\, \le \, {\mathrm{Tr}}\, (f(A)+f(B)).$$
}

\vskip 5pt\noindent In the matrix setting, the concavity assumption
is quite crucial as shown in the following simple remark \cite{Ty}:

\vskip 10pt\noindent
 {\it Let
$f:[0,\infty)\to[0,\infty)$ be continuous with $f(0)=0$. If
$${\mathrm{Tr}}\, f(A+B)\, \le \, {\mathrm{Tr}}\, (f(A)+f(B))$$
holds for all two-by-two positive matrices $A,\,B$, then $f$ is
concave.}

\vskip 10pt\noindent To prove this statement, take for $s,t>0$,
$$
A=\frac{1}{2}
\begin{pmatrix} s&\sqrt{st} \\
\sqrt{st}&t
\end{pmatrix}
\qquad B=\frac{1}{2}
\begin{pmatrix} s&-\sqrt{st} \\
-\sqrt{st}&t
\end{pmatrix}
$$
and observe that the trace inequality means that $f$ is concave.

 Theorem 1.1  closed a list of papers of several authors including
 Ando-Zhan \cite{AZ}, and Kosem \cite{Ko}.
However, It remained natural to ask wether this result could be
extended to the set of all Hermitian, or even  all normal operators.
We noticed  a partial answer in \cite{B3}:

\vskip 10pt\noindent
 {\bf Theorem 1.2.} {\it Let $A,\,B\ge0$ and let
 $f:[0,\infty)\to [0,\infty)$ be concave and e-convex. Then, for all symmetric norms,
 $$\Vert\, f(|A+B|)\,\Vert \le \Vert\, f(|A|)+f(|B|)\,\Vert.$$
 }

\vskip 5pt \noindent Here the  e-convexity property of $f$ means
that $f(e^t)$ is convex on $(-\infty,\infty)$. In particular, the
theorem holds for the power functions $f(t)=t^p$, $1\ge p\ge 0$.
This result for normal operators   entails several  estimates for
block matrices. A special case involving an operator partitioned in
four normal blocks $A,\, B,\, C,\, D$ of same size  is
\begin{equation*}
\left\| \,\left|
\begin{pmatrix}A&B \\C&D \end{pmatrix}
 \right|^p \,\right\|
  \le \left\|\, |A|^p+|B|^p+|C|^p+|D|^p\,\right\|
\end{equation*}
for all symmetric norms and $0\le p\le 1$.  These estimates,
comparing an operator on ${\mathcal{H}}\oplus {\mathcal{H}}$ with a
related operator on ${\mathcal{H}}$, differ from the usual ones in
the literature where the norm of the full matrix is evaluated with
the norms of its blocks, for instance, \cite{KgNa} and \cite{BhKi1}.
In the subsequent sections we solve the conjectures in \cite{B3} by
showing that the assumption of e-convexity is not necessary in
Theorem 1.2,   in its application to block-matrices and in some
related inequalities. The proof of Theorem 1.2 given in [7] reduced
to the positive case by using the fact that for any normal $A$, $B$
and any non-negative e-convex functions $f(t)$, we have
$$
\Vert\, f(|A+B|)\,\Vert \le \Vert\, f(|A|+|B|)\,\Vert
$$
for all symmetric norms. This is no longer true if the e-convexity
assumption is dropped. In fact, one can easily find two-by-two
positive semi-definite matrices $A$, $B$ and a non-negative concave
function $f(t)$ on $[0,\infty)$ such that
\begin{equation*}
\Vert\, f(|A-B|)\,\Vert > \Vert\, f(|A|+|B|)\,\Vert
\end{equation*}
for all symmetric norms which are not a multiple scalar of the usual
operator norm. For instance, take $f(t)=\min\{t,\sqrt2/2\}$ and
\begin{equation}
A=\begin{pmatrix}1&0 \\0&0 \end{pmatrix}, \qquad B=
\begin{pmatrix}1/2&1/2 \\1/2&1/2 \end{pmatrix}.
\end{equation}
The main point of the forthcoming proof is to overcome this
difficulty. This proof can be adapted in order to obtain a version
for normal operators of the following companion result to Theorem
1.1:

\vskip 10pt\noindent
 {\bf Theorem 1.3.} {\it Let $A\ge0$ and let  $Z$  be
 expansive.
 If
 $f:[0,\infty)\to [0,\infty)$  is concave, then, for all symmetric norms,
 $$\Vert\, f(Z^*AZ)\,\Vert \le \Vert\, Z^*f(A)Z\,\Vert.$$
}

 \vskip 15pt\noindent {\Large\bf 2. Subadditivity results for normal operators}

\vskip 15pt We have the following norm inequalities:

\vskip 10pt\noindent
 {\bf Theorem 2.1.} {\it Let $A,\,B$ be normal and let
 $f:[0,\infty)\to [0,\infty)$ be concave. Then, for all symmetric norms,
 $$\Vert\, f(|A+B|)\,\Vert \le \Vert\, f(|A|)+f(|B|)\,\Vert.$$
 }

\vskip 10pt\noindent
 {\bf Corollary 2.2.}  {\it Let   $Z=A+iB$ be a  decomposition in real and imaginary parts, and let
 $f:[0,\infty)\to [0,\infty)$ be concave. Then, for all symmetric norms,
$$
\left\|\,f(|Z|)\, \right\| \le \left\|\,f(|A|)+f(|B|)\,\right\|.
$$
}

\vskip 5pt\noindent
 This is a matrix version of the obvious
inequality $f(|z|)\le f(|a|)+f(|b|)$ for  complex numbers $z=a+ib$.
Since non-negative concave functions on $[0,\infty)$ are
non-decreasing we actually have $f(|z|)\le f(|a|+|b|)\le
f(|a|)+f(|b|)$. But the left inequality can not be extended to
matrices. Indeed it is easy to find two-by-two matrices $Z=A+iB$ - a
simple example is given with $A$, $B$ defined in (2) - with the
eigenvalue relation
$$
\lambda_2(|A|+|B|)<\lambda_2(|Z|)<\lambda_1(|Z|)<\lambda_1(|A|+|B|).
$$
Thus, there are some non-negative concave functions like
$f(t)=\min\{t,\lambda_2(|Z|)\}$ such that
$$
\Vert\, f(|Z|)\,\Vert > \Vert\, f(|A|+|B|)\,\Vert
$$
for all symmetric norms which are not a multiple scalar of the usual
operator norm.

Let $A$, $B$ be general operators. Applying Theorem 2.1 to the
Hermitian operators
$$\begin{pmatrix} 0&A^*\\ A&0\end{pmatrix} \qquad \mathrm{and}
\qquad \begin{pmatrix} 0&B^*\\ B&0\end{pmatrix}
$$
we obtain
$$
\left\|\begin{pmatrix} f(|A+B|)&0\\
0&f(|A^*+B^*|)\end{pmatrix}\right\| \le
\left\|\begin{pmatrix} f(|A|)+f(|B|)&0\\
0&f(|A^*|)+f(|B^*|)\end{pmatrix}\right\|
$$
so that,
letting $B= A^*$ yields:

\vskip 10pt\noindent
 {\bf Corollary 2.3.}  {\it If
 $f:[0,\infty)\to [0,\infty)$ is concave,
 then, for all $Z$ and all symmetric norms,
$$
\Vert \,f(|Z+Z^*|)\,\Vert \le  \Vert \,f(|Z|)+f(|Z^*|)\,\Vert.
$$
} \vskip 5pt\noindent Note that equality occurs in Corollary 2.3
whenever $f(0)=0$ and
$$
Z=\begin{pmatrix} 0&0\\ X&0\end{pmatrix}
$$
where $X$ is arbitrary. Note also that it may happen that
$$
\Vert \,f(|Z+Z^*|)\,\Vert >  \Vert \,f(|Z|+|Z^*|)\,\Vert
$$
for some concave functions and some symmetric norms, for instance
when
$$
Z=\begin{pmatrix}0&1&0 \\0&0&1 \\ 0&0&0 \end{pmatrix}
$$
and the norm is the sum of the two largest singular values.

At the end of this section, we will see some application of Theorem
2.1 to partitioned operators. Now, we turn to the proof of Theorem
2.1. We start by recalling the Ky Fan Principle. The Ky Fan
$k$-norms of $A$, $k=1,\,2,\ldots,n$ are defined as the sum of its
$k$ largest singular values,
$$
\Vert A\Vert_{(k)}=\sum_{j=1}^k\lambda_j(|A|).
$$
Let $A$, $B$ such that $\Vert A\Vert_{(k)}\le \Vert B\Vert_{(k)}$
for all $k=1,\,2,\ldots,n$. Then, the vector of the singular values
of $A$ lies in the convex hull of the permuted singular values of
$B$ multiplied by $\pm 1$,
$$
(\lambda_1(|A|),\cdots,\lambda_n(|A|)) \in {\mathrm{conv}}_{\sigma}
(\pm \lambda_{\sigma(1)}(| B|),\cdots,\pm \lambda_{\sigma(n)}(|B|))
$$
This can be proved by using the  Hyperplan separation process to
reach a contradiction, see \cite{Si} for details and \cite{Bh} for
alternative proofs. From this convexity statement follows a useful
fact:

\vskip 10pt\noindent
 {\bf Ky Fan Principle.} {\it Suppose that $\Vert A\Vert_{(k)}\le \Vert B\Vert_{(k)}$ for
 all Ky-Fan $k$-norms. Then, we  have $\Vert A\Vert \le \Vert B\Vert$ for all symmetric
 norms.
 }

\vskip 10pt\noindent
 We also need two elementary, well-known lemmas. For $A,\,B\ge 0$ it is sometimes
 convenient to write $A\prec_w B$ to mean that $\Vert A\Vert \le \Vert B\Vert$ for all symmetric
 norms.

\vskip 10pt\noindent
 {\bf Lemma 1.} {\it Let $A,\,B,\,X,\,Y \ge 0$ such that $B\prec_w Y$
 and $A\prec_w X$. Then,
$$
\begin{pmatrix}A&0 \\0&B \end{pmatrix}\prec_w \begin{pmatrix}X&0 \\0&Y
\end{pmatrix}.
$$
}

\vskip 10pt\noindent
 {\bf Proof.} We have
 $$
 \sum_{j=1}^k\lambda_j(A\oplus B)=\max_{s+t=k}\left\{\sum_{j=1}^s\lambda_j(A) +
 \sum_{j=1}^t\lambda_j(B)\right\}.
 $$
  Combining this with
 $$
 \sum_{j=1}^s\lambda_j(A) + \sum_{j=1}^t\lambda_j(B)\le  \sum_{j=1}^s\lambda_j(X) + \sum_{j=1}^t\lambda_j(Y)\le  \sum_{j=1}^k\lambda_j(X\oplus Y)
 $$
 ends the proof. \qquad $\Box$

\vskip 10pt\noindent
 {\bf Lemma 2.} {\it Let $A,\,B \ge 0$. Then,
$$
\begin{pmatrix}A&0 \\0&B \end{pmatrix}\prec_w \begin{pmatrix}A+B&0
\\0&0
\end{pmatrix}.
$$
}

\vskip 10pt\noindent {\bf Proof.} Note that
$$
\begin{pmatrix}A+B&0\\0&0
\end{pmatrix}=
\begin{pmatrix}A^{1/2}&B^{1/2}\\0&0
\end{pmatrix}
\begin{pmatrix}A^{1/2}&0\\B^{1/2}&0
\end{pmatrix}
$$
so that
$$
\begin{pmatrix}A+B&0\\0&0
\end{pmatrix}
\simeq
\begin{pmatrix}A&A^{1/2}B^{1/2}\\B^{1/2}A^{1/2}&B
\end{pmatrix}
\simeq
\begin{pmatrix}A&-A^{1/2}B^{1/2}\\-B^{1/2}A^{1/2}&B
\end{pmatrix}
$$
where $\simeq$ means unitarily congruent. Combining with
$$
\begin{pmatrix}A&0 \\0&B \end{pmatrix}
=\frac{1}{2}\begin{pmatrix}A&A^{1/2}B^{1/2}\\B^{1/2}A^{1/2}&B\end{pmatrix}
+\frac{1}{2}\begin{pmatrix}A&-A^{1/2}B^{1/2}\\-B^{1/2}A^{1/2}&B\end{pmatrix}
$$
gives the lemma. \qquad $\Box$

 \vskip 10pt\noindent
 {\bf Proof of Theorem 2.1.}  It  suffices to prove the result when $A$ and $B$ are Hermitian. The general case then follows by
 replacing $A$, $B$ by
$$
 \tilde{A}=\begin{pmatrix}0&A \\A^*&0 \end{pmatrix}, \qquad \tilde{B}=\begin{pmatrix}0&A \\A^*&0 \end{pmatrix}
$$
and by using normality of  $A$ and $B$. Therefore assume that $A$,
$B$ are Hermitian with decomposition in positive and negative parts,
$$A=A_+-A_- \qquad {\mathrm{and}}\qquad  B=B_+-B_-.$$
Let $g(t)=f(t)-f(0)$ and note that, for each Ky Fan $k$-norm,
$$
\Vert\, f(|A+B|)\,\Vert_{(k)} = kf(0)+\Vert\, g(|A+B|)\,\Vert_{(k)}
$$
and
$$
\Vert\, f(|A|)+f(|B|)\,\Vert_{(k)} = 2kf(0)+\Vert\,
g(|A|)+g(|B|)\,\Vert_{(k)}.
$$
Hence, it suffices to prove the result for $g(t)$, or equivalently
when $f(0)=0$. This assumption implies
\begin{equation}
f(|A|)=f(A_+)+f(A_-)\qquad {\mathrm{and}} \qquad
f(|B|)=f(B_+)+f(B_-).
\end{equation}
 Now, given two positive $n$-by-$n$ matrices
$X$ and $Y$ with direct sum
$$
X\oplus Y=\begin{pmatrix} X&0\\ 0&Y\end{pmatrix}
$$
 we have
\begin{equation}
\lambda_j(|X-Y|)\le \lambda_j(X\oplus Y)
\end{equation}
for all $j=1,\,2,\ldots, n$. Indeed, for some subspace
${\mathcal{S}}\subset{\mathcal{H}}$ we have
\begin{align*}
|X-Y|&=(X-Y)_+ +(X-Y)_-\\
&=(X-Y)_{{\mathcal{S}}}\oplus(Y-X)_{{\mathcal{S}}^{\perp}}\\
&\le X_{{\mathcal{S}}}\oplus Y_{{\mathcal{S}}^{\perp}}
\end{align*}
hence
$$
\lambda_j(|X-Y|)\le \lambda_j(X_{{\mathcal{S}}}\oplus
Y_{{\mathcal{S}}^{\perp}})\le\lambda_j(X\oplus Y)
$$
 for all $j=1,\,2,\ldots, n$. Replacing in (4) $X$ by $A_{+}+B_{+}$
and $Y$ by $A_-+B_-$ we then get
\begin{equation*}
\lambda_j(|A+B|)\le \lambda_j((A_{+}+B_{+})\oplus (A_-+B_-))
\end{equation*}
for all $j=1,\,2,\ldots, n$. Since $f$ is non-decreasing, it follows
\begin{equation*}
\lambda_j(f(|A+B|))\le \lambda_j(f(A_{+}+B_{+}))\oplus (f(A_-+B_-))
\end{equation*}
for all $j=1,\,2,\ldots, n$, so that
$$
\Vert\,f(|A+B|)\,\Vert \le \left\|\begin{pmatrix} f(A_{+}+B_{+}) &0\\
0&f(A_-+B_-)
\end{pmatrix}\right\|
$$
for all symmetric norms. By Theorem 1.1 combined with Lemma 1,
followed by application of Lemma 2, we then obtain
$$
\Vert\,f(|A+B|) \,\Vert\le
\Vert\,f(A_{+})+f(B_{+})+f(A_{-})+f(B_{-})\,\Vert
$$
and making use of relations (3) ends the proof. \qquad $\Box$

 \vskip 10pt
 Let us now give some application for Block-matrices. The most
 obvious one is for a Hermitian matrix
\begin{equation*}
 \begin{pmatrix}A&B \\B^*&C \end{pmatrix}
 \end{equation*}
partitioned in four blocks of same size. Then by using Theorem 2.1
for the decomposition in two Hermitian
\begin{equation*}
 \begin{pmatrix}A&B \\B^*&C \end{pmatrix}=
\begin{pmatrix}A&0 \\ 0&C \end{pmatrix}+\begin{pmatrix}0&B \\B^*&0 \end{pmatrix}
 \end{equation*}
  and then using Lemma 1, we have
\begin{equation}
 \left\| \,f\left(\left|
\begin{pmatrix}A&B \\B^*&C \end{pmatrix}\right|
  \right)\,\right\|
  \le \left\|\, f(|A|)+f(|B|)+f(|B^*|)+f(|C|)\,\right\|
\end{equation}
for all concave functions $f:[0,\infty)\to [0,\infty)$ and all
symmetric norms.

To obtain similar statements for more general partitions, note that
the  proof of Theorem 2.1 is valid for any finite family of normal
operators. Thus: {\it  Let $\{A_i\}_{i=1}^m$ be normal and let
$f:[0,\infty)\to [0,\infty)$ be concave. Then, for all symmetric
norms,
  $$\Vert\, f(|A_1+\cdots+A_m|)\,\Vert \le \Vert\, f(|A_1|)+\cdots+f(|A_m|)\,\Vert.$$
 }
 We may then obtain results for some matrices partitioned in $m^2$
blocks of same size.

 \vskip 10pt\noindent
 {\bf Corollary 2.4.}  {\it Let ${\mathbb{A}}=[A_{i,\,j}]$ be a block matrix with normal
 entries and let  $f$ be a non-negative concave function on $[0,\infty)$.
  Then, for all symmetric norms,
$$
\left\|\, f(|{\mathbb{A}}|) \,\right\| \le \left\|\,\sum
f(|A_{i,\,j}|) \,\right\|.
$$
}

 \vskip 10pt\noindent {\bf Proof.}  We prove this corollary  via Theorem 2.1. for a partition in four
blocks
$$
{\mathbb{A}}=\begin{pmatrix}S&R \\T&Q \end{pmatrix}.
$$
The proof for a   partition in $m^2$ blocks is similar by using the
version of Theorem 2.1 for $m$ operators. Let
$$
\tilde{{\mathbb{A}}}=\begin{pmatrix}0&{\mathbb{A}}
\\{\mathbb{A}}^*&0
\end{pmatrix}
$$
and note that
\begin{equation*}
|\tilde{{\mathbb{A}}}|=\begin{pmatrix}|{\mathbb{A}}^*|&0
\\ 0&|{\mathbb{A}}|\end{pmatrix}
\end{equation*}
so that
\begin{equation}
|\tilde{{\mathbb{A}}}|\simeq\begin{pmatrix}|{\mathbb{A}}|&0
\\ 0&|{\mathbb{A}}|\end{pmatrix}
\end{equation}
where the symbol $\simeq$ stands for unitarily equivalent. On the
other hand
 $$
 \tilde{{\mathbb{A}}}=\tilde{{\mathbb{S}}}+\tilde{{\mathbb{T}}}
 $$
 where
 $$
\tilde{{\mathbb{S}}}=\begin{pmatrix} 0&0&S&0 \\ 0&0&0&Q
\\ S^*&0&0&0 \\ 0&Q^*&0&0\end{pmatrix}
\quad \tilde{{\mathbb{T}}}=\begin{pmatrix} 0&0&0&R \\ 0&0&T&0
\\ 0&T^*&0&0 \\ R^*&0&0&0\end{pmatrix}
$$
 are Hermitian. Therefore, Theorem 2.1 yields,
 $$
\left\|\,f(|\tilde{{\mathbb{A}}}|)\,\right\|\,\le\,
\left\|\,f(|\tilde{{\mathbb{S}}}|)+f(|\tilde{{\mathbb{T}}}|)
\,\right\|
$$
for all symmetric norms; that is, using the shorthand symbol
$\prec_w$,
$$f(|\tilde{{\mathbb{A}}}|)\prec_w
 \begin{pmatrix} f(|S^*|)+f(|R^*|)&0&0&0 \\ 0&f(|T^*|)+f(|Q^*|)&0&0
\\ 0&0&f(|S|)+f(|T|)&0 \\ 0&0&0&f(|R|)+f(|Q|)\end{pmatrix}.
$$
 Gathering the two first lines, and the two last ones, we have via
Lemmas 2 and 1
$$
 f(|\tilde{{\mathbb{A}}}|)\prec_{w}
 \begin{pmatrix} f(|S^*|)+f(|T^*|)+f(|R^*|)+f(|Q^*|)&0\\ 0&f(|S|)+f(|T|)+f(|R|)+f(|Q|)\end{pmatrix}.
$$
By using (6) we then obtain, using normality of $S,\,T,\,R,\,Q$,
$$
f(|{\mathbb{A}}|)\prec_{w} f(|S|)+f(|T|)+f(|R|)+f(|Q|)
$$
which is equivalent to inequalities for symmetric norms. \qquad
$\Box$

\vskip 10pt Let us point out a variation of Corollary 2.4 in which
some operators are not necessarily normal.

\vskip 10pt\noindent
 {\bf Corollary 2.5.} {\it Let ${\mathbb{T}}$ be a triangular block-matrix
 $$
 {\mathbb{T}}=\begin{pmatrix} A&N\\0&B\end{pmatrix}.
 $$
 in which $N$ is normal. Let $f:[0,\infty)\to[0,\infty)$ be  concave. Then, for all symmetric
 norms,
 $$
\Vert f(|{\mathbb{T}}|)\Vert \le \Vert f(|A^*|)+ f(|N|)+f(|B|)
\Vert.
 $$
 }

\vskip 10pt\noindent
 {\bf Proof.} Consider the polar decompositions $A=|A^*|U$ and
 $B=V|B|$, note that
 $$
\left|\begin{pmatrix} A&N\\0&B\end{pmatrix}\right|\simeq \left|
\begin{pmatrix} I&0\\0&V^*\end{pmatrix}\begin{pmatrix}
A&N\\0&B\end{pmatrix}\begin{pmatrix}
U^*&0\\0&I\end{pmatrix}\right|=\left|\begin{pmatrix}
|A^*|&N\\0&|B|\end{pmatrix}\right|
 $$
 and apply Theorem 2.1. \qquad $\Box$

 \vskip 10pt The
assumption in Corollary 2.4 requiring normality of each block is
rather special. The next corollary generalizes (5) and meets the
simple requirement that the full matrix is Hermitian.

 \vskip 15pt\noindent
 {\bf Corollary 2.6.}  {\it Let ${\mathbb{A}}=[A_{i,\,j}]$ be a  Hermitian matrix  partitioned in blocks of same size and let
  $f:[0,\infty)\to[0,\infty)$ be  concave.
  Then, for all symmetric norms,
$$
\left\|\, f(|{\mathbb{A}}|) \,\right\| \le \left\|\,\sum
f(|A_{i,\,j}|) \,\right\|.
$$
}

\vskip 10pt\noindent {\bf Proof.} The proof of Corollary 2.4
actually shows that for a general block-matrix
${\mathbb{A}}=(A_{i,\,j})$ partitioned in blocks of same size, we
have
$$
 \begin{pmatrix}f(|{\mathbb{A}}|)&0\\ 0&f(|{\mathbb{A}}|)\end{pmatrix}\prec_{w}
 \begin{pmatrix} \sum f(|A_{i,\,j}^*|)&0\\ 0& \sum f(|A_{i,\,j}|)\end{pmatrix}
$$
for all non-negative concave function $f$.
Assuming ${\mathbb{A}}$ Hermitian, we have $A_{i,\,j}^*=A_{j,\,i}$
and Corollary 2.6 follows.  \qquad $\Box$

\vskip 25pt\noindent {\Large\bf 3. Related results for expansive
congruences}

\vskip 15pt Let $A$ be normal and let $Z$ be expansive, i.e.,
$Z^*Z\ge I$. The following extension of Theorem 1.3 holds.

\vskip 10pt\noindent
 {\bf Theorem 3.1.} {\it Let $A$ be normal and let  $Z$  be
 expansive.  If
 $f:[0,\infty)\to [0,\infty) $  is concave, then, for all symmetric norms,
 $$\Vert\, f(|Z^*AZ|)\,\Vert \le \Vert\, Z^*f(|A|)Z\,\Vert.$$
}

\vskip 10pt\noindent Indeed, we can derive Theorem 3.1 from Theorem
1.3 in a quite similar way of the one we derive Theorem 2.1 from
Theorem 1.1. The proof  of Theorem 3.1 starts by noticing that we
can assume that $A$ is Hermitian. Then, using the decomposition in
positive and negative parts
 $$A=A_+-A_- $$
 we have, as in the proof of Theorem 2.1,
 $$
 \lambda_j(|Z^*AZ|) \le\lambda_j(|Z^*A_+Z|\oplus |Z^*A_-Z|)
 $$
 and we may proceed as previously.

 When we deal with the trace norm, the fact that $f$ is positive on
 the whole half-line is not essential, as in Rotfel'd inequality.
 Hence we have the following corollary, extending to normal operators
 a result from \cite{B1}.

\vskip 10pt\noindent {\bf Corollary 3.2.} {\it Let $A$ be normal and
let $Z$ be expansive. If $f(t)$ is a concave function on the
positive half-line with $f(0)\ge0$, then
 $${\mathrm{Tr}}\, f(|Z^*AZ|)\, \le \,{\mathrm{Tr}}\, Z^*f(|A|)Z.$$
}

\vskip  5pt\noindent Corollary 3.2 follows from Theorem 3.1 by
approaching $f(t)$ with $g(t)+at$ for some scalar $a$ and some
non-negative concave function $g(t)$. Theorem 2.1 and 3.1 can be
combined in a unique statement, extending the main result in
\cite{BL}:

 \vskip 10pt\noindent
 {\bf Theorem 3.3.}  {\it Let $\{A_i\}_{i=1}^m$ be normal, let  $\{Z_i\}_{i=1}^m$ be expansive and let
 $f$ be a non-negative concave function on $[0,\infty)$. Then, for all symmetric norms,
$$
\left\|f\left(\left|\sum Z_i^*A_iZ_i\right|\right) \right\| \le
\left\| \sum Z_i^*f(|A_i|)Z_i \right\|.
$$
}

\vskip 10pt\noindent It would be elegant and interesting to state
this theorem in the more general  framework of positive linear maps
$\Phi$ between matrix algebras. This leads to the problem of
characterizing the positive linear maps $\Phi$ such that
 \begin{equation*}
 \Vert f(|\Phi(N)|)\Vert \le \Vert \Phi( f(|N|))\Vert
 \end{equation*}
for all normal operators $N$, all non-negative concave functions and
all symmetric norms. Some furher questions are considered in
\cite{B3}. For sake of completeness, we mention that when $f(t)$ is
a non-negative convex function vanishing at 0, then inequalities of
Theorem 1.1-1.3 are reversed. For instance we have \cite{BL}

\vskip 10pt\noindent
 {\bf Theorem 3.4.}  {\it Let $\{A_i\}_{i=1}^m$ be positive and  let  $\{Z_i\}_{i=1}^m$ be expansive.
  Then, for all symmetric norms and all $p>1$,
$$
\left\| \sum Z_i^*A_i^pZ_i \right\| \le \left\|  \left(\sum
Z_i^*A_iZ_i \right)^p \right\|.
$$
}

\vskip 5pt\noindent
 If $Z_i=I$ for all $i$, it is a famous result
of Ando-Zhan [1] and of Bhatia-Kittaneh [3] in case of integer
exponents. The very special case $ {\mathrm Tr\,} (A_1^p+A_2^p) \le
{\mathrm Tr\,} (A_1+A_2) ^p$ is Mc-Carthy's inequality [13, p.\ 20].
Note that the  positivity assumption in Theorem 3.4 can not be
replaced by a normality one.

When we consider contractive congruences and positive operators,
then there exist several Jensen type inequalities, not only for
norms but also for eigenvalues (cf.\ \cite{B1} \cite{B2}). The proof
are much simpler than in the expansive case, where some unexpected
counterexamples may occur (see discussion and counterexamples in
\cite{B1} \cite{B2}). We give an example of such results:

\vskip 10pt\noindent
 {\it Let $\{A_i\}_{i=1}^m$ be positive and
$\{Z_i\}_{i=1}^m$ such that $\sum Z^*_iZ_i \le I$. If
 $f$ is a monotone concave function on $[0,\infty)$, $f(0)\ge0$, then,
\begin{equation*}
f\left(\sum Z_i^*A_iZ_i\right)  \ge   V\left(\sum Z_i^*f(A_i)Z_i
\right)V^*.
\end{equation*}
 for some unitary $V$.}

 \vskip 15pt

%%%%%%%%%%%%%%%%%%%%%%%%%%%%%%%%%%%%%%%%%%%%%%%%%%%%%%%%%%%%%%%%%%%%%%%%%%%

\vskip 20pt\noindent

 Jean-Christophe Bourin

 jcbourin@@univ-fcomte.fr

 Laboratoire de math\'ematiques

 Universit\'e de Franche-Comt\'e

 25030 Besan\c{c}on

\end{document}